\documentclass[twocolumn,showpacs,preprintnumbers,amsmath,amssymb]{revtex4}

\usepackage{psfig} 

\usepackage{graphicx} 

\begin{document}

\preprint{APS/123-QED}

\title{Improving Search Algorithms by Using Intelligent Coordinates}

\author{David Wolpert}
\author{Kagan Tumer}
\author{Esfandiar Bandari}%
\affiliation{NASA Ames Research Center, Moffett Field, CA, 94035, USA}

\begin{abstract}
We consider the problem of designing a set of computational agents so
that as they all pursue their self-interests a global function $G$ of
the collective system is optimized. Three factors govern the quality
of such design.  The first relates to conventional
exploration-exploitation search algorithms for finding the maxima of
such a global function, e.g., simulated annealing (SA). Game-theoretic
algorithms instead are related to the second of those factors, and the
third is related to techniques from the field of machine learning.
Here we demonstrate how to exploit all three factors by modifying the
search algorithm's exploration stage so that rather than by random
sampling, each coordinate of the underlying search space is controlled
by an associated machine-learning-based ``player'' engaged in a
non-cooperative game.  Experiments demonstrate that this modification
improves SA by up to an order of magnitude for bin-packing and for a
model of an economic process run over an underlying network. These
experiments also reveal novel small worlds
phenomena.
\end{abstract}



\maketitle              

\section{INTRODUCTION}
\label{sec:intro}

Many distributed systems found in nature have inspired
function-maximization algorithms. In some of these the coordinates of
the underlying system are viewed as players engaged in a
non-cooperative game, whose joint behavior (hopefully) maximizes the
pre-specified global function of the entire system. Examples of such
systems are auctions and clearing of markets. Typically in the
computer-based algorithms inspired by such ``collectives'' of players,
each separate coordinate of the system is controlled by an associated
machine learning
algorithm \cite{cama97,chjo02,crba96,huho88,sama97,zhan98},
reinforcement-learning (RL) algorithms being particularly
common~\cite{suba98,wowh00}.

One important issue concerning such collectives is whether the payoff
function $g_\eta$ of each player $\eta$ is sufficiently sensitive to
what coordinates $\eta$ controls in comparison to the other
coordinates, so that $\eta$ can learn how to control its coordinates
to achieve high payoff.  A second crucial issue is the need for all of
the $g_\eta$ to be ``aligned'' with $G$, so that as the players
individually learn how to increase their payoffs, $G$ also increases.

Previous work in the COllective INtelligence (COIN) framework
addresses these issues. This work extends conventional game-theoretic
mechanism design~\cite{futi91, niro01} to include off-equilibrium
behavior, learnability issues, $g_\eta$ with non-human attributes
(e.g., $g_\eta$ for which incentive compatibility is irrelevant), and
arbitrary $G$.  In domains from network routing to congestion problems
it outperform traditional techniques, by up to several orders of
magnitude for large
systems~\cite{tuwo00,wowh00,wotu01a}.

Other collective systems found in nature that have inspired search
algorithms do not involve players conducting a non-cooperative game.
Examples include spin glasses, genomes undergoing neo-Darwinian
natural selection, and eusocial insect colonies, which have been
translated into simulated annealing (SA ~\cite{kige83,gege84}),
genetic algorithms~\cite{bafo97,chfo99}, and swarm
intelligence~\cite{bodo00,krbi00}, respectively.  An important issue
here is the exploration/exploitation dynamics of the overall
collective.

Recent analysis reveals how $G$ is governed by the interaction between
exploration/exploitation, the alignment of the $g_\eta$ and $G$, and
the learnability of the $g_\eta$ ~\cite{wotu01a}.  Here we use that
analysis to motivate a hybrid algorithm, Intelligent Coordinates for
search (IC), that addresses all three issues.  It works by modifying
any exploration-based search algorithm so that each coordinate being
searched is made ``intelligent'', its exploration value being the move
of a game-playing computer algorithm rather than the random sample of
a probability distribution.

Like SA, IC is intended to be used as an ``off the shelf'' algorithm;
rarely will it be the best possible algorithm for some particular
domain.  Rather it is designed for use in very large problems where
parallelization can provide a large advantage, while there is little
exploitable information concerning gradients.  We present experiments
comparing IC and SA on two archetypal domains: bin-packing and an
economic model of people choosing formats for their home music
systems.

In the bin-packing domain IC achieves a given value of $G$ up to three
orders of magnitude faster than does SA, with the improvement
increasing linearly with the size of the problem.  In the format
choice problem $G$ is the sum of each person's ``happiness'' with her
format choices. Each person $\eta$'s happiness with each of her
choices is set by three factors: which of her nearest neighbors on a
ring network ($\eta$'s ``friends'') make that choice; $\eta$'s
intrinsic preference for that choice; and the price of music purchased
in that format, inversely proportional to the total number of players
using that choice. Here again, IC improves $G$ two orders of magnitude
more quickly than does SA.  We also considered an algorithm similar to
the Groves mechanism of economics; IC outperformed it by over two
orders of magnitude.  We also modified the ring to be a small-worlds
network ~\cite{kual00,nemo00,wast98}.  This barely improved IC's
performance (3\%), with no effect on the other algorithms.  However if
$G$ was also changed, so that each $\eta$'s happiness depends on
agreeing with her friends' friends, the performance increase in
changing to a small-worlds topology is significant (10\%). This
underscores the multiplicity of factors behind the benefits of
small-worlds networks.

\section{Simplified Theory of Collectives}
\label{sec:math}
Let $z \in \zeta$ be the joint move of all agents/players in the
collective.  We want the $z$ that maximizes the provided world utility
$G(z)$. In addition to $G$ we have private utility functions
\{$g_\eta$\}, one for each agent $\eta$ controlling $z_\eta$.
$\hat{}\;\eta$ refers to all agents other than $\eta$.

{\bf Intelligence} ``standardizes'' utility functions so that the
value they assign to $z$ only reflects their ranking of $z$ relative
to some other $z'$. One form of it is
\begin{eqnarray}
N_{\eta,U}(z ) \equiv \int d\mu_{z_{\;\hat{}\;\eta}} (z')
	\Theta[U(z) - U(z')] \; ,
\end{eqnarray}

\noindent
where $\Theta$ is the Heaviside function, and where the subscript on
the (normalized) measure $d\mu$ indicates it is restricted to $z'$
such that $z'_{\hat{}\;\eta} = z_{\hat{}\;\eta}$.

Our uncertainty concerning the system induces a distribution
$P(z)$. All attributes of the collective we can set, e.g., the private
utility functions of the agents, are given by the value of the {\bf
design coordinate} $s$. Bayes theorem provides the {\bf central
equation}:

\begin{eqnarray}
& & P(G \mid s) = \\ \nonumber
& &  \int  d\vec{N}_G P(G \mid \vec{N}_G, s) \int
d\vec{N}_g P(\vec{N}_G \mid \vec{N}_g, s) P(\vec{N}_g \mid s) \; ,
\label{eq:central}
\end{eqnarray}

\noindent
where $\vec{N}_G$ and $\vec{N}_g$ are the intelligence vectors for all
the agents, for utilities $g_\eta$ and for $G$, respectively.
$N_{\eta,g_\eta}(z) = 1$ means that agent $\eta$'s move maximizes its
utility, given the moves of the other agents.  So $\vec{N}_{g}(z) =
\vec{1}$ means $z$ is a Nash equilibrium.  Conversely, $\vec{N}_G(z')
= \vec{1}$ means that the value of $G$ cannot increase in moving from
$z'$ along any single (sic) coordinate of $\zeta$. So if these two
points are identical, then if the agents do well enough at maximizing
their private utilities they must be near an (on-axis) maximizing
point for $G$.

More formally, say for our $s$ the third conditional probability in
the integrand in the central equation (``term 3'') is peaked near
$\vec{N}_g = \vec{1}$. Then $s$ probably induces large (private
utility function) intelligences (intuitively, the utilities are
learnable).  If in addition the second term is peaked near $\vec{N}_G
= \vec{N}_g$, then $\vec{N}_G$ will also be large (intuitively, the
private utility is ``aligned with $G$''). This peakedness is assured
if $\vec{N}_g = \vec{N}_G$ exactly $\forall z$. Such a system is said
to be {\bf factored}.  Finally, if the first term in the integrand is
peaked about high $G$ when $\vec{N}_G$ is large, then $s$ probably
induces high $G$, as desired.

As a trivial example, a {\bf team game}, where $g_\eta = G \; \forall
\eta$, is factored~\cite{crba96}.  However team games usually have
poor third terms, especially in large collectives.  This is because
each $\eta$ has to discern how its moves affect $g_\eta = G$, given
the background of the (varying) moves of the other agents whose moves
comparably affect $G$.

Fix some $f(z_\eta)$, two moves ${z_\eta}^1$ and ${z_\eta}^2$, a
utility $U$, a value $s$, and a $z_{\;\hat{}\;\eta}$. The associated
{\bf learnability} is
\begin{eqnarray}
\label{eq:learnability}
& & \Lambda_f(U ; z_{\;\hat{}\;\eta}, s, {z_\eta}^1, {z_\eta}^2)
\equiv \\ \nonumber
& & \sqrt{ \frac{ [E(U ; z_{\;\hat{}\;\eta}, {z_\eta}^1) -
 E(U ; z_{\;\hat{}\;\eta}, {z_\eta}^2)]^2}
{\int dz_\eta [f(z_\eta)Var(U ; z_{\;\hat{}\;\eta}, z_\eta)] }} \;
 .
\end{eqnarray}
\noindent
The averages and variance here are evaluated according to $P(U |
n_\eta) P(n_\eta | z_{\;\hat{}\;\eta}, {z_\eta}^1)$, $P(U | n_\eta)
P(n_\eta | z_{\;\hat{}\;\eta}, {z_\eta})$, and $P(U | n_\eta) P(n_\eta
| z_{\;\hat{}\;\eta}, {z_\eta}^2)$, respectively, where $n_\eta$ is
$\eta$'s training set, formed by sampling $U$.

The denominator in Eq.~\ref{eq:learnability} reflects the sensitivity
of $U(z)$ to $z_{\;\hat{}\;\eta}$, while the numerator reflects its
sensitivity is to $z_{\eta}$. So the greater the learnability of
$g_\eta$, the more $g_\eta(z)$ depends only on the move of agent
$\eta$, i.e., the more learnable $g_\eta$ is.  More formally, it can
be shown that if appropriately scaled, $g'_\eta$ will result in better
expected intelligence for agent $\eta$ than will $g_\eta$ whenever
$\Lambda_f(g'_\eta; z_{\;\hat{}\;\eta}, s, {z_\eta}^1, {z_\eta}^2) \;
> \; \Lambda_f(g_\eta; z_{\;\hat{}\;\eta}, s, {z_\eta}^1, {z_\eta}^2)$
for all pairs of moves ${z_\eta}^1, {z_\eta}^2$\cite{wolp03a}.

A {\bf difference} utility is one of the form $U(z) = G(z) -
D(z_{\;\hat{}\;\eta})$.  Any difference utility is
factored ~\cite{wolp03a}. In addition, under usually benign
approximations, the $D(z_{\;\hat{}\;\eta})$ that maximizes
$\Lambda_f(U ; z_{\;\hat{}\;\eta}, s, {z_\eta}^1, {z_\eta}^2)$ for all
pairs ${z_\eta}^1, {z_\eta}^2$ is $E_f(G(z) \mid z_{\;\hat{}\;\eta},
s)$, where the expectation value is over $z_\eta$. The associated
difference utility is called the {\bf Aristocrat} utility ($AU$).  If
each $\eta$ uses its $AU$ as its private utility, then we have both
good terms 2 and 3.

Evaluating the expectation value in $AU$ can be difficult in practice.
This motivates the {\bf Wonderful Life} Utility (WLU), which requires
no such evaluation:
\begin{eqnarray}
WLU_\eta \equiv G(z) - G(z_{\;\hat{}\;\eta}, CL_\eta) \; ,
\end{eqnarray}
\noindent
where $CL_\eta$ is the {\bf clamping parameter}.  WLU is factored,
independent of the clamping parameter. Furthermore, while not matching
$AU$, $WLU$ typically has far better learnability than does a team
game, and therefore typically results in better values of $G$. It is
also often easier to evaluate than is $G$ itself \cite{tuwo00,wotu01a}.

One way to address term 1 as well as 2 and 3 is to incorporate
exploration/exploitation techniques like SA.

\section{EXPERIMENTS}

In our version of SA, at the beginning of each time-step $t$ a
distribution $h_\eta(\zeta_{\eta})$ is formed for every $\eta$ by
allotting probability 75\% to the move $\eta$ had at the end of the
preceding time-step, $z_{\eta,t-1}$, and uniformly dividing
probability 25\% across all of its other moves. The ``exploration''
joint-move $z_{expl}$ is then formed by simultaneously sampling all
the $h_\eta$. If $G(z_{expl}) > G(z_{t-1})$, $z_{\eta,t}$ is set to
$z_{expl}$. Otherwise $z_{t}$ is set by sampling a Boltzmann
distribution having energies $G(z_{t-1})$ and $G(z_{expl})$. Many
different annealing schedules were investigated; all results below are
for best schedules found.

IC is identical except that each $h_\eta$ is replaced by
$\frac{h_\eta(z_\eta)c_{\eta,t}(z_\eta)} {\sum_{a'}
h_\eta(a'_\eta)c_{\eta,t}(a'_\eta)}$, where the distribution
$c_{\eta,t}$ is set by an RL algorithm trying to optimize payoffs
$g_\eta$.  Here RL is done using a training set $n_{\eta,t}$ of all
preceding move-payoff pairs, \{$(z_{\eta,t'}, g_\eta(z_{t'}) : t' <
t)$\}. For each possible move by $\eta$ one forms the weighted average
of the payoffs recorded in $n_{\eta,t}$ that occurred with that move,
where the weights decay exponentially in $t - t'$.  $c_{\eta,t}$ then
is the Boltzmann distribution, parameterized by a ``learning
temperature'' (that effectively rescales $g_\eta$) over those
averages.

In all our experiments the ``AU'' version of IC approximated $f$ to to
be uniform $\forall \eta$, and then used a mean-field approxmation to
pull the expectation inside $G$.  Unless otherwise specified, the
clamping elements used in WLU's were set to $\vec{0}$.



In bin-packing $N$ items, all of size $< c$, must be assigned into a
minimal subset of $N$ bins, without assigning a summed size $> c$ to
any one bin.  $G$ of an assignment pattern is the number of occupied
bins~\cite{coga98}, and each agent controls the bin choice of one
item.  To improve performance all algorithms use a modified ``$G$'',
$G_{\mbox{soft}}$, even though their performance is measured with $G$:
\begin{eqnarray}
G_{\mbox{soft}} = 
\left\{ \begin{array}{ll}
      \sum_{i=1}^N   \left[ \left(\frac{c}{2}\right)^2  - 
      \left(  x_i - \frac{c}{2}\right)^2 \right] 
					  & \mbox{if $x_i \leq c$}\\
      \sum_{i=1}^N
      \left(  x_i  - \frac{c}{2}\right)^2 & \mbox{if $x_i > c$}
      \end{array}
\right.
,
\label{eq:Gbinpack}
\end{eqnarray}
where $x_i$ is the summed size of all items in bin $i$.  (Use of
$G_{\mbox{soft}}$ encourages bins to be either full or empty.)

In the IC runs learning temperature was .2, and all agents made the
transition to RL-based moves after a period of 100 random $z$'s used
to generate the starting $n_{\eta}$.  Exploitation temperature started
at .5 for all algorithms, and was multiplied by .8 every 100
exploitation time-steps In each SA run, the distribution $h$ was
slowly modified to generate solutions that differed in fewer items
than the current solution as time progressed.

\begin{table}[htb]  \centering
\begin{tabular}{c|c|c|c|c} \hline
Algorithm  & Ave. $G$     &  Best & Worst & \% Optimum
\\ \hline
IC WLU   & 3.32 $\pm$ 0.22 & 2 &  8  & 72 \% \\
IC TG    & 7.84 $\pm$ 0.17 & 6 & 10  & 0  \% \\ \hline
COIN WLU & 3.52 $\pm$ 0.20 & 2 &  7  & 64 \% \\
COIN TG  & 7.84 $\pm$ 0.15 & 6 &  9  & 0  \% \\ \hline
SA       & 6.00 $\pm$ 0.19 & 4 &  7  & 0  \% \\ \hline
\end{tabular}
\caption{Bin-packing $G$ at time 200 for $N=20,c=12$.}
\label{tab:binpack}
\end{table}

\begin{figure}[bth]
\vspace*{-.2in}
\centerline{\psfig{figure={ 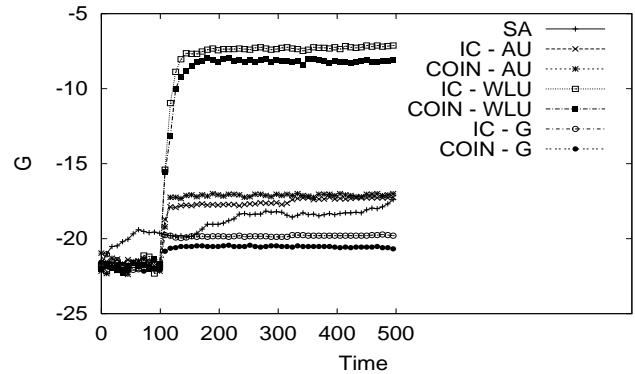},width=3.5in,height=2.2in}}
\caption{Average bin-packing $G$ for $N=50, c=10$. All error bars
$\leq .31$ except IC - AU and COIN - AU are $\leq .57$.}
\label{fig:bin50}
\end{figure}

\begin{figure}[bth]
\vspace*{-.2in}
\centerline{\psfig{figure={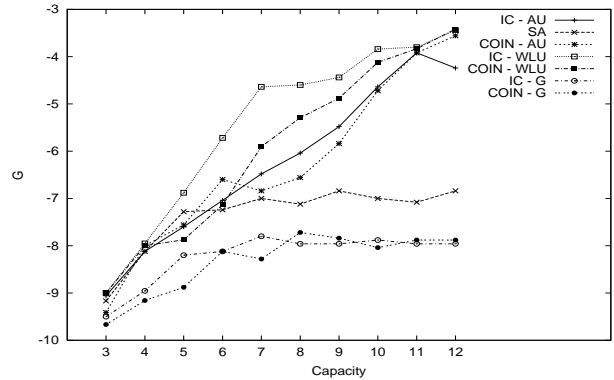},width=3.3in,height=2.0in}}
\caption{$G$ vs. $c$ for $N=20$ at $t=200$. All
error bars $\leq .34$.}
\label{fig:cap}
\end{figure}

In Table 1 ``Best'' refers to the best end-of-run $G$ achieved by the
associated algorithm, ``worst'' is the worst value, and ``\%Optimum''
is the percentage of runs that were within one bin of the best value.
Fig.~\ref{fig:bin50} shows average performances (over 25 runs) as a
function of time step. The algorithms that account for both terms 2
and 3 --- IC WLU and COIN WLU --- far outperform the others, with the
algorithm accounting for all three terms doing best.  The worst
algorithms were those that accounted for only a single term (SA and
COIN TG).  Linearly (i.e., optimistically) extrapolating SA's
performance from time 15000 indicates it would take over 1000 times as
long as IC WLU to reach the $G$ value IC WLU reaches at time 200. In
addition the ratio of WLU's time 1000 performance (relative to random
search) to SA's grows linearly with the size of the problem. Finally,
Fig.~\ref{fig:cap} illustrates that the benefit of addressing terms 2
and 3 grows with the difficulty of the problem. In both
figures SA outperforms IC - TG; this is due to there being more
parameter-tuning with SA.

For the format choice problem $G$ is the sum over all $N_a$ agents
$\eta$ of $\eta$'s ``happiness'' with its music formats:
\begin{equation}
G = \sum_{\eta=1}^{N_a}
\sum_{i=1}^{N_f} \sum_{\eta' \in neigh_\eta} \vartheta(i) 
\; \omega_{\eta,\eta',i} \; pref_{\eta,i}
\end{equation}
where $N_f$ is the numbers of formats; $neigh_\eta$ is the set of
players lying $\le D$ hops away from player $\eta$; $pref_{\eta,i}$ is
$\eta$'s intrinsic preferenece for format $i$ (set randomly at
initialization $\in [0, 1]$); $\vartheta_(i)$ is the total number of
players that choose format $i$ (i.e., the inverse price for format
$i$); and $\omega_{i,\eta,\eta'} = 1$ if the choices of players $\eta$
and $\eta'$ both include the format $i$, and 0 otherwise (each agent's
move is a selection of three of four total formats, implemented by
choosing the one format not to be used). $D$ values of both 1 and 3
were investigated.


\begin{figure}[bth]
\centerline{\psfig{figure={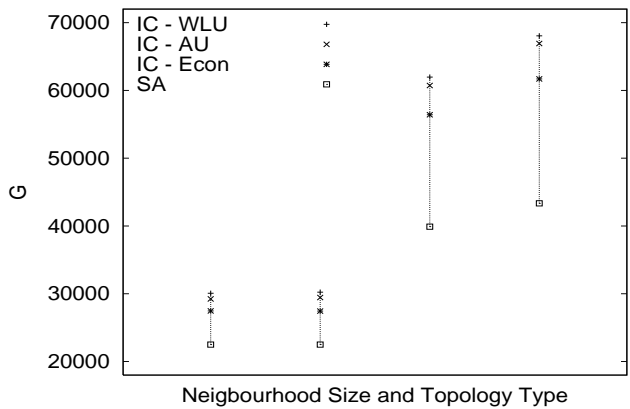},width=3.5in,height=2.2in}}


\caption{$G(t = 200)$ for 100 agents. In order from left to right, $D
= \{1,1,3,3\}$, and topologies are \{L,W,L,W\}.}
\label{fig:sw100}
\end{figure}

In Fig.~\ref{fig:sw100}, ``IC Econ'' refers to WLU IC where clamping
means the agent chooses no format whatsoever. It is essentially the
game-theory Groves mechanism wherein one sets $g_\eta$ to $\eta$'s
marginal contribution to $G$, here rescaled and interleaved with a
simulated annealing step to improve performance. ``IC-WLU'' instead
clamps $\eta$'s move to zero (in accord with the theory of
collectives), which means that $\eta$ chooses all formats.  Learning
temperature was now .4, and exploitation temperature was .05
(annealing provided no advantage since runs were short).  Two network
topologies were investigated. Both were $m$-node rings with an extra
$.06 m$ random links added, a new such set for each of the 50 runs
giving a plotted average value. ``Short links'' (L) means that all
extra links connected players two hops apart, while ``small-worlds''
(W) means there was no such restriction.


IC Econ's inferior performance illustrates the shortcoming of
economics-like algorithms. For $D=1$ SA did not benefit from small
worlds connections, and IC variants barely benefited (~ 3\%), despite
the associated drop in average inter-node hop distance.  However if
$D$ also increased, so that $G$ directly reflected the change in the
topology, then the gain with a small worlds topology grew to ~
10\%. (See the discussion on path lengths in~\cite{nemo00}.)





The authors thank Michael New, Bill Macready, and Charlie Strauss for
helpful comments.

\bibliographystyle{plain}

\end{document}